\documentclass[12pt]{amsart}
\usepackage{amssymb}
\usepackage{amsmath}
\newtheorem{thm}{Theorem}[section]
\newtheorem{lemma}[thm]{Lemma}
\newtheorem{theorem}[thm]{Theorem}
\newtheorem*{corollary*}{Corollary}

\theoremstyle{definition}
\newtheorem{definition}[thm]{Definition}

\newtheorem{remark}[thm]{Remark}

\numberwithin{equation}{section}

\newcommand{\CC}{{\mathbb C}}

\newcommand{\bbC}{{\mathbb C}}

\newcommand{\fg}{{\mathfrak g}}

\renewcommand{\div}{\mbox{div}}
\renewcommand{\Im}{\mbox{Im}}

\newcommand{\eff}{\text{eff}}
\newcommand{\cut}{\text{cut}}

\renewcommand{\Re}{\mbox{Re}}   
\newcommand{\Span}{\mbox{span}}
\newcommand{\stable}{\text{stable}}

\newcommand{\grad}{\nabla}

\def \bar{\overline}
\def \hat{\widehat}

\numberwithin{equation}{section}

\theoremstyle{definition}
\theoremstyle{remark}

\newcommand{\R}{\mathbb{R}}
\newcommand{\C}{\mathbb{C}}
\newcommand{\Z}{\mathbb{Z}}

\newcommand{\bbT}{\mathbb{T}}

\def	\inv	{^{-1}}

\newcommand{\Diff}{\mathop{\it Diff}\nolimits}

\newcommand\bM {\bar{M}} 

\begin{document}

\title{K\"ahler cuts}

\author{D. Burns \and V. Guillemin \and E. Lerman}\thanks{Supported in
part by NSF grants DMS-0104047 (DB), DMS-0104116 (VG) and
DMS-0204448(EL)} 


\address{M.I.T., Cambridge, MA 02139 \and University of Michigan, Ann
Arbor, MI 48109}

\email{dburns@umich.edu}

\address{M.I.T., Cambridge, MA 02139}

\email{vwg@math.mit.edu}

\address{University of Illinois, Urbana, IL 61801}

\email{lerman@math.uiuc.edu}

\begin{abstract}

A symplectic cut of a manifold M with a Hamiltonian circle action is a
symplectic quotient of $M \times \mathbb{C}$.  If $M$ is K\"ahler then,
since $\mathbb{C}$ is K\"ahler, the cut space is K\"ahler as well.
The symplectic structure on the cut is well understood.  In this paper
we describe the complex structure (and hence the metric) on the cut.
We then generalize the construction to the case where $M$ has a torus
action and $\mathbb{C}$ is replaced by a toric K\"ahler manifold.

\end{abstract}

\maketitle

\section{Introduction}
\label{sec:1.}

Let $(M,\omega)$ be a symplectic manifold and $\tau: S^1\times M\to M$
a Hamiltonian action of $S^1$ with moment map $\phi : M \to
\mathbb{R}$. We will assume that, for $\lambda \in \mathbb{R}$, $S^1$
acts freely on $\phi^{-1}(\lambda)$. In particular, $\lambda$ is then
a regular value of $\phi$ and the symplectic quotient
\begin{equation}
M_{\lambda}=\phi^{-1}(\lambda)/S^1 
\end{equation}
is well-defined and non-singular. In \cite{Le} it was shown that the
disjoint union of this quotient with the open subset 
\begin{equation}
\label{eqn:M^lambda_o}
M^\lambda_o = \{ p \in M \mid \phi(p) < \lambda\} 
\end{equation}
of $M$ can be given the structure of a smooth symplectic manifold
$(M^\lambda, \omega^\lambda)$ which was called the \emph {symplectic
cut} of $M$ at $\lambda$. A number of applications of this cutting
operation were also given to problems in symplectic geometry, and
since then many other applications have been found. (See, for
instance, \cite{DGMW}, \cite{LMTW}, \cite{PW}, \cite{MS}, and
\cite{MW}.)

In this article we will look at this cutting operation from the
K\"ahlerian perspective. Our motivation is a basic problem in K\"ahler
geometry: what happens to the K\"ahler metric on a non-singular
projective variety when one blows this variety up along a non-singular
subvariety? The symplectic cutting operation turns out to have some
bearing on this problem. Suppose that the moment map above takes its
maximum value $\lambda_o$ on a subset $W$ of $M$. Then $W$ is a
symplectic submanifold of $M$, and, for $\lambda = \lambda_o - \epsilon,
\epsilon \approx 0$, $M^\lambda$ can be obtained from $M$ by blowing
up $M$ symplectically along $W$. (See \cite{Mc} and \cite{GS}.)

Suppose now that $M$ is a complex manifold and $\omega$ a K\"ahler
form.  In addition, suppose that $\tau$ extends to a holomorphic action
of $\mathbb{C}^*$ on $M$. Then $M^\lambda$ is also K\"ahler, so the set
(\ref{eqn:M^lambda_o}) has two K\"ahler structures: the K\"ahler
structure it acquires as an open subset of $M^\lambda$ and the
K\"ahler structure it acquires as an open subset of $M$. In \cite{Le}
it was shown that symplectically these two structures coincide:
$\omega = \omega^\lambda$ on $M^\lambda_o$. However, it was also
pointed out that the two complex structures on $M^\lambda_o$
\emph{don't} coincide. The main result of this paper is an
amplification of this remark, a global description of the complex
structure on $M^\lambda_o$ coming from $M^\lambda$. Let 
\begin{equation}
\label{eqn:Msharp}
M^\# = \mathbb{C}^*\cdot M^{\lambda}_o. 
\end{equation}
For $M$ compact this is a Zariski open subset of $M$, and we will
construct below a canonical diffeomorphism of $M^\lambda_o$ onto
$M^\#$ which is biholomorphic with respect to the $M^\lambda$ complex
structure on $M^\lambda_o$ and the $M$ complex structure on $M^\#$.

Since $M \sp \#$ is Zariski open, its complement is a complex
subvariety of $M$, and we will show that this variety is the union of
unstable manifolds for the $\mathbb{C} \sp{*}$-action, or, in other
words, ``generalized Schubert varieties". Thus the cutting operation
compactifies the complement of these Schubert varieties by adjoining
the non-singular hypersurface $M \sb \lambda$ to $M$.

We give now a brief summary of the contents of this article. In
$\S~2$ we review the definition and elementary properties of the
cutting operation and prove the assertions above. In $\S~3$ we
generalize the version of the cutting operation given in  \cite{Le} by
replacing the standard K\"ahler potential $|z|^2$ on
$\C$ by an arbitrary radial potential $F$, and show that, with minor
modifications, the results of $\S~2$ are still true. Then, in $\S~4$ we
apply these results to the K\"ahler potential $F(z) = \frac{2}{c}
\log (1 + |z|^2), c > 0$ and show that if the K\"ahler structure
on $M$ is K\"ahler-Einstein with structure constant $c$, the K\"ahler
form on $M \sp \lambda$ satisfies a modified version of the
K\"ahler-Einstein equation. Finally, in $\S~5$ we explain how to
generalize the results of this paper to the ``toric variety" version
of symplectic cutting developed in \cite{Le} and \cite{LMTW}.

\section{Symplectic cuts}
\label{sec:2.}

Let $(M,\omega)$ be a K\"ahler manifold and $\tau:S^1\times M \to M$ a
Hamiltonian action of $S^1$ on $M$ with moment map $\phi:M \to
\R$.  We will assume that this action extends to a holomorphic
action of $\C^*$ on $M$ which we will continue to denote by
$\tau$. Suppose that $S^1$ acts freely on the level
set $\phi^{-1}(\lambda)$. Let $W = M \times\C$, equipped with
the K\"ahler form 
\begin{equation}
\label{eqn:omegacut}
\omega + \sqrt{-1}dz\wedge d\bar{z}
\end{equation}
and let $S^1$ act on $W$ by its product action. This is a Hamiltonian
action with moment map 
\begin{equation}
\label{eqn:psicut}
\psi = \phi + |z|^2 .
\end{equation}
Moreover, $S^1$ acts freely on the level set, $\psi^{-1}(\lambda)$, so the
reduced space 
\begin{equation}
\label{eqn:cutdef}
M^\lambda = \psi^{-1}(\lambda)/S^1
\end{equation}
is well-defined, and this symplectic manifold is by definition the
manifold ``$M$ cut at $\lambda$''(see \cite{Le}). The level set
$\psi^{-1}({\lambda})$. is the disjoint union of the two $S^1$-invariant
sets: 
\begin{equation}
\label{eqn:philevel}
\{(p,0) \mid \phi(p) = \lambda\}
\end{equation}
and
\begin{equation}
\label{eqn:phisublevel}
\{(p,z) \mid \phi(p) = \lambda -|z|^2 < \lambda\},
\end{equation}
so the quotient $\psi^{-1}(\lambda)/S^1$ is the disjoint union of
the quotients of these two sets. The quotient of \ref{eqn:philevel} is
by definition the reduced space
$$
M_\lambda = \phi^{-1}(\lambda)/S^1
$$
and the quotient of (\ref{eqn:phisublevel}) can be identified with the
open subset 
\begin{equation}
M^\lambda_o = \{p\in M \mid \phi(p) < \lambda\}
\end{equation}
of $M$. In fact this identification is given explicitly by the map
\begin{equation}
\label{eqn:map}
M^\lambda_o \ni p \to (p, \sqrt{\lambda - \phi(p)})
\end{equation}
which maps $M^\lambda_o$ onto a global cross-section for the action of
$S^1$ on the set (\ref{eqn:phisublevel}). Thus $M^\lambda$ is the
disjoint union of the 
open subset $M^\lambda_o$ of $M$ and a codimension two symplectic
submanifold, the reduced space $M_\lambda$. In particular,
$M^\lambda_o$ has two K\"ahler structures: the restriction of the
K\"ahler structure on $M^\lambda$ and the restriction of the K\"ahler
structure on $M$. We propose to show below how these
K\"ahler structures are related. 

>From the symplectic point of view the
situation is extremely simple. Let $\omega^\lambda$ be the symplectic
form on $M^\lambda$. The pull-back by the map (\ref{eqn:map}) of the
form (\ref{eqn:omegacut}) is just $\omega$, and so on $M^\lambda_o$,
$\omega = \omega^\lambda$. However, as is pointed out in \cite{Le}, section
1.1, $(M^\lambda_o,\omega)$ and $(M^\lambda_o, \omega^\lambda)$ are
\emph{not} identical as K\"ahler manifolds: the $M$-complex structure
on $M^\lambda_o$ doesn't coincide with the $M^\lambda$-complex
structure. How these complex structures are related is the question we
will investigate below. 

Let $v = \grad \phi$ be the infinitesimal generator of the one parameter 
group $\tau_{e^t}:M\to M$, and for $q\in M^\lambda_o$ let
\begin{equation}
\label{eqn:kappadef}
\kappa(q) = -1/2 \log(\lambda - \phi(q)).
\end{equation}
We define a map $g:M^\lambda_o \to M$ by setting 
\begin{equation}
\label{eqn:mapg}
g(q) = \exp(\kappa(q) v)(q).
\end{equation}
Let $M_{\stable} = \mathbb{C}^*\cdot\phi^{-1}(\lambda)$. We claim that
$g$ is a diffeomorphism of $M^\lambda_o$ onto the open subset
\begin{equation}
\label{eqn:Msharpbis}
M^\# = M^\lambda_o \cup M_{\stable}
\end{equation}
of $M$. That $g$ is a diffeomorphism onto its image is clear. To see
that its image is (2.8) we first note that if $q$ is in $M^\lambda_o$,
then $p = g(q)$ if and only if 
\begin{equation}
\label{eqn:veqn}
\exp(tv)(p) = q 
\end{equation}
for $t$ satisfying 
\begin{equation}
\label{eqn:teqn}
\lambda - e^{2t} = \phi(q).
\end{equation}

Now let $p$ be in $M_{\stable}$ with $\phi(p) > \lambda$, and let
$\gamma(t) = \exp(tv)(p)$. Then
\begin{equation}
a = \lim_{t \to -\infty}\phi(\gamma(t)) < \lambda ,
\end{equation}
so the curves, $y = \lambda - e^{2t}$ and $y = \phi(\gamma(t))$ must
intersect at some point on the negative $t$ axis; and, at that
point $q =\gamma(t)$ satisfies $p = g(q)$ by (\ref{eqn:veqn}) and
(\ref{eqn:teqn}). Similarly if $p$ is in $M^\lambda_o$ and $\gamma(t) =
\exp(tv)(p)$, then if $\phi(p) < \lambda - 1$, the curve, $y =
\phi(\gamma(t))$, has to intersect the curve,$y = \lambda - e^{2t}$ at
some point on the positive $t$ axis, and if $\lambda - 1 < \phi(p)
<\lambda$, these curves must intersect at some point on the
negative $t$ axis. Thus the image of $g$ contains $M^\lambda_o \cup
M_{\stable}$, and the inclusion the other way is obvious.  If $M$ is
compact the set (\ref{eqn:Msharpbis}) doesn't depend on $\lambda$ but only
on the critical values of $\phi$ lying above $\lambda$. More
explicitly, let $F_i, i = 1,...,k$, be the connected components of the
set of critical points of $\phi$, i.e., the connected components of
the fixed point set of the action $\tau$. Let $W_i$ be the unstable
manifold of $\grad \phi$ at $F_i$.  Thus
\begin{equation}
q\in W_i \Leftrightarrow \lim_{z \to 0} \tau_z(q)
\in F_i 
\end{equation}
These manifolds are complex submanifolds of positive codimension. Let
$S$ be the set of $i$'s for which $\phi(F_i)$ is greater than
$\lambda$. We claim
\begin{equation}
\label{eqn:decomposition}
M -M^\# = \bigcup_{i\in S} W_i.
\end{equation}
\noindent Proof: By the Whitney decomposition theorem $M$ is the
disjoint union of the $W_i$'s so every point $p \in M$ is in some
$W_i$, and if $p$ is in $W_i$ and $\phi(F_i)$ is less than $\lambda$,
then either $\phi(p) <\lambda$ or the $\mathbb{C}^*$ orbit through $p$
intersects $\phi^{-1} (\lambda)$. We will now prove:

\begin{thm}\label{thm2.1}
The map $g: M^\lambda_o \to M^\#$ is a biholomorphism of $M^\lambda_o$ with its
$M^\lambda$ complex structure onto $M^\#$ with its $M$ complex structure.
\end{thm}

\begin{proof}

The biholomorphic map
\begin{equation}
\label{eqn:mapf}
f: M\times \mathbb{C}^* \to  M\times \mathbb{C}^*
\end{equation}
defined by
\begin{equation}
\label{eqn:fformula}
f(p,z_o) = (\tau_{z_o}p,z_o)
\end{equation}
intertwines the action
\begin{equation}
\label{eqn:ract}
z(p,z_o) = (p,zz_o)
\end{equation}
of $\mathbb{C}^*$ on $M\times\mathbb{C}^*$ with the diagonal action
\begin{equation}
\label{eqn:dact}
z(p,z_o) = (\tau_zp,zz_o)
\end{equation}
so the pullback by $f$ of the K\"ahler form (2.1) is invariant under the
action of $S^1$ on the second factor of $M\times\mathbb{C}^*$. Moreover
this action is Hamiltonian with moment map
\begin{equation}
\label{eqn:dacthamiltonian}
\tilde{\psi}(p,z) = \phi(\tau_zp) + |z|^2
\end{equation}
In particular, $f$ maps the level set, $\tilde{\psi}^{-1}(\lambda)$
onto the level set $\psi^{-1}(\lambda)$, and induces an isomorphism of
K\"ahler manifolds
\begin{equation}
\label{eqn:maph}
h:\tilde\psi^{-1}(\lambda)/S^1 \to \psi^{-1}(\lambda)/S^1 
\end{equation}
To describe this isomorphism more explicitly note that the set
\begin{equation}
\{\Re \, z > 0, \Im \, z = 0\}
\end{equation}
is a global cross-section for both of the $S^1$ actions above, so
$\tilde\psi^{-1}(\lambda)/S^1$ can be identified with the set
\begin{equation}
\label{eqn:ractsect}
\{(p,e^t) \mid \phi(\tau_{e^t}p) + e^{2t} = \lambda\}
\end{equation}
and $\psi^{-1}(\lambda)$ with the set
\begin{equation}
\label{eqn:dactsect}
\{(q,e^t) \mid \phi(q) + e^{2t} =\lambda\}
\end{equation}
and, modulo these identifications, $h$ maps the point $(p,e^t)$ in
the set (\ref{eqn:ractsect}) onto the point $(q,e^t)$ in the set
(\ref{eqn:dactsect}) where 
$$
q= \tau_{e^t}p
$$
and
$$
\lambda - e^{2t} = \phi(q).
$$
However, $\tau_{e^t}p = \exp(tv)(p)$, so by (\ref{eqn:teqn}) $p = g(q)$,
i.e., $h$ is just $g^{-1}$. In particular, the domain of $h$ is the open
subset $M^\#$ of $M$ defined by (\ref{eqn:Msharpbis}). Identifying $M^\#$
with $\tilde\psi^{-1}(\lambda)/S^1$, we note that since $\mathbb{C}^*$ is
acting on $M \times \mathbb{C}^*$ by the action (\ref{eqn:dact}) the
geometric invariant theory quotient 
$$
M \times \mathbb{C}^*/\!/\mathbb{C}^* = \tilde\psi^{-1}(\lambda)/S^1
$$
is not just equal to the open set $M^\#$ set theoretically, but is
this open set with its $M$-complex structure. (See, for instance,
\cite{BG}, \S~4.) Thus $g$ intertwines the $M^\lambda$-complex
structure on $M^\lambda_o$ with the complex structure on $M$.

\end{proof}

\noindent Example: Let $\lambda_o$ be the maximum value of $\phi$ and let
$W = \phi^{-1}(\lambda_o)$. Then $W$ is a complex submanifold of $M$, and
if $\lambda = \lambda_o - \epsilon, \epsilon \approx 0$, then 
$$
M^\# = M - W
$$
by \ref{eqn:decomposition}. On the other hand, it is not hard to show
that $M_\lambda$ is the projective normal bundle, $\mathbb{P}(NW)$,
of $W$. Therefore, the assertion that $M^\lambda -M_\lambda$ is
biholomorphic to $M^\#$ implies that $M^\lambda$ can be obtained from
$M$ by deleting $W$ and adjoining $\mathbb{P}(NW)$, i.e., by blowing up $M$ along $W$.

Next we will describe some equivariance properties of the mapping $g$. 
The holomorphic action of $\mathbb{C}^*$ on $M\times \mathbb{C}^*$
defined by (\ref{eqn:dact}) extends to a holomorphic product action of
$\mathbb{C}^*\times\mathbb{C}^*$ 
\begin{equation}
\label{eqn:pact}
(z_1,z_2)\cdot(p,z) = (\tau_{z_1}p, z_2z)
\end{equation}
and from this one gets a residual holomorphic action of $\mathbb{C}^*$
on the quotient space $M^\lambda$. The submanifold $M_\lambda$ of
$M^\lambda$ is fixed by this action, and hence one gets an induced
action of $\mathbb{C}^*$ on its complement $M^\lambda_o$. Restricted
to $S^1$ this action coincides with the $M$-action of $S^1$ on
$M^\lambda_o$; however, this action itself can't coincide with the
$M$-action of $\mathbb{C}^*$ since the $M$-action of $\mathbb{C}^*$
doesn't leave $M^\lambda_o$ fixed.\footnote{This is one of the more
compelling ways of seeing that the $M$-complex structure on
$M^\lambda_o$ can't be identical with the $M^\lambda$-complex
structure.} As we noted above, however, the action of $\mathbb{C}^*$
on $M$ does leave fixed the submanifold $M^\#$. We claim

\begin{thm}

The mapping $g$ intertwines the $M$-action of $\mathbb{C}^*$ on
$M^\#$ with the $M^\lambda$ action of $\mathbb{C}^*$ on $M^\lambda_o$

\end{thm}

\begin{proof} 

Since the function (2.7) and the vector field $v$ are $S^1$-invariant,
the map $g$ intertwines the two $S^1$ actions. Therefore, since $g$ is
biholomorphic, it intertwines the two $\mathbb{C}^*$ actions as
well. 

\end{proof}

Let $U$ be a $\mathbb{C}^*$-invariant open subset of $M$, and let
$(z_1, ...,z_n)$ be a holomorphic coordinate system on $U$. We will
say that say that $(z_1, ...,z_n)$ is a $\mathbb{C}^*$-\emph{adapted
coordinate system} if $z_n \neq 0$ on $U$ and if the action of $\mathbb{C}^*$ on $U$ is given by
\begin{equation}
\tau_a(z_1, ...,z_n) = (z_1, ...,z_{n-1}, az_n).
\end{equation}
Then, for $p\in M^\lambda_o \cap U$
\begin{equation}
g(p) = \tau_{e^t}(p), t = \kappa(p).
\end{equation}
But $\kappa(p) = -1/2 \log (\lambda -\phi(p))$ so
\begin{equation}
\label{eqn:coor}
p = (z_1, ...,z_n) \Leftrightarrow g(p) = 
(z_1, ...,z_{n-1},(\lambda - \phi(z))^{-1/2}z_n).
\end{equation}
In other words
\begin{equation}
g^*z_i = z_i, i = 1, ..., n-1,
\end{equation}
and
$$
g^*z_n = (\lambda - \phi(z))^{-1/2}z_n
$$
Thus from theorem 2.1 we conclude

\begin{thm}

If $U$ is a $\mathbb{C}^*$-invariant open subset of $M$ and $(z_1,
...,z_n)$ a $\mathbb{C}^*$-adapted coordinate system on $U$, then
(\ref{eqn:coor}) is a complex coordinate system on $U \cap M^\lambda_o$
compatible with the $M^\lambda$-complex structure.
\end{thm}

\section{Cuts with invariant potentials}
\label{sec:3.}

In the construction of the cut space $M^\lambda$ described in $\S~2$
above one began by equiping $W = M \times \CC$ with the K\"ahler form
the sum of that on $M$ and the form $\sqrt{-1} dz \wedge d\bar{z}$ on
$\CC$, and then reduced with respect to the diagonal action of
$S^1$. There is, of course, nothing sacrosanct about the form
$\sqrt{-1} dz \wedge d\bar{z}$, and in fact, it is sometimes
advantageous to consider more general K\"ahler forms on the $\CC$
factor, {\em e.g.}, when $M$ is K\"ahler-Einstein, ({\em cf.}, $\S~4$). In
this section we will replace $\sqrt{-1} dz \wedge d\bar{z}$ by an
{\em arbitrary} $S^1$-invariant K\"ahler form on $\CC$ and show that
most of the results of $\S~2$ are still true, modulo small changes
detailed below.

Let $s$ be a real valued function on $\C$ invariant under the standard
action of $S^1$: $\lambda \cdot z = \lambda z$.  By a theorem of
G. Schwarz there exists a function $F$ on $\R$ such that $s(z,\bar{z})
= F (|z|^2)$.  Now 
\begin{equation}
\partial {\bar \partial} F (|z|^2) =
\left( F'' (|z|^2) |z|^2 + F' (|z|^2) \right) dz \wedge d\bar{z}.
\end{equation}
Let $H(t) = tF(t)$.
\begin{lemma}
Let $F$ and $H$ be as above. Then
\begin{enumerate}
\item  $\omega_F:= \sqrt{-1}\partial {\bar \partial} F (|z|^2)$ is a 
symplectic form on $\C$ iff $H'(t) >0$ for all $t\geq 0$.

\item $\phi = H (|z|^2)$ is a moment map for the standard action of
$S^1$ on $(\C, \omega_F)$.

\end{enumerate}
\end{lemma}

\begin{proof}
The form $\omega_F$ is nondegenerate iff $F'' (|z|^2) |z|^2 + F'
(|z|^2) > 0 $ for all $z\in \C$ iff $ H'(t) = F'' (t) t + F'(t) > 0$
for all $t\geq 0$.   This proves the first claim.
Our proof of the second claim is a computation:
$$
\iota \left( \sqrt{-1} (z\frac{\partial}{\partial z} - \bar{z}\frac{\partial}
{\partial \bar{z}})\right) \omega_F = 
- ( F'' (|z|^2) |z|^2 + F' (|z|^2)) (z d\bar{z} + \bar{z} dz) =
 - d H(|z|^2).
$$
\end{proof}
 
Note that by definition $H(0) = 0$.  Also, since $H'(t) > 0$ on $[0,
\infty)$, $H(t)$ is strictly increasing on $[0, \infty)$.  Let $a =
\lim_{t\to + \infty} H(t)$, so that $0 < a \leq \infty$.  Then $H$ is
invertible and its inverse $K$ is a strictly increasing function $K :
[0, a) \to [0, \infty)$.

Now let $(M, \omega)$ be a symplectic manifold with a Hamiltonian
action $\tau$ of $S^1$ and let $\phi: M \to \R$ be an associated
moment map. For $\lambda \in \R$ we define the {\bf cut}  of $(M,
\omega)$ with respect to the potential $F$ at $\lambda$  to be the symplectic
quotient at $\lambda$ of $(M\times \CC, \omega + \omega _F)$ under the
diagonal action of $S^1$.  We denote the cut by $M_{\text{cut}} (
\lambda, F)$. Of course $M^\lambda$ (\ref{eqn:cutdef}) of \S~2 above
is just $M_{\cut}(\lambda, |z|^2)$.

\begin{thm}\label{thm3.2}

Let $(M, \omega)$ be a symplectic manifold with a Hamiltonian action
of $S^1$ and let $\phi: M \to \R$ be an associated moment map.
Suppose $S^1$ acts freely on $\phi \inv (\lambda)$.  Then the cut
$M_{\cut} ( \lambda, F)$ is naturally a smooth symplectic manifold.

Moreover, the natural embedding of the reduced space $M_{\lambda} :=
\phi\inv (\lambda)/S^1$ into $M_{\cut} ( \lambda, F)$ is
symplectic, and the complement $M_{\cut} ( \lambda, F)
- M_{\lambda}$ is symplectomorphic to the open subset $M^{\lambda}_o :=
\{m\in M \mid \lambda -a < \phi (m)< \lambda\}$ of $(M, \omega)$, where as above
$a = \lim_{|z|\to + \infty} \phi(z)$.

Furthermore, if $(M, \omega)$ is a K\"ahler manifold and the action of
$S^1$ is holomorphic, then $M_{\text{cut}} ( \lambda, F)$ is K\"ahler
and $M_{\text{cut}} ( \lambda, F) - M_\lambda$ is
biholomorphic to $M^\# := \C^* \cdot M^{\lambda}_o = \C^* \cdot \{m\in
M \mid \lambda -a < \phi(m)< \lambda\} \subset M$.    

\end{thm} 

\begin{proof}

The moment map $\psi$ for the diagonal action of $S^1$ on $(M \times
\C, \omega + \omega_F)$  is given by 
$\psi (m ,z) = \phi (m) + H (|z|^2)$ and so
\begin{equation}
\label{eqn:psilevelbis}
\psi \inv (\lambda) = 
\{ (m, z) \in M\times \C \mid \phi (m) + H (|z|^2) = \lambda\}.
\end{equation}
Now $\phi (m) + H (|z|^2) = \lambda$ iff $|z|^2 = K ( \lambda -
\phi(m))$, where as before $K= H\inv$.  Consider the map
\begin{equation}
\label{eqn:sigmadef}
\sigma :  \{ m \in M \mid \lambda -a < \phi(m)\leq \lambda\} \to \psi
\inv (\lambda), \; \sigma (m) = (m, (K (\lambda -\phi
(m)))^{\frac{1}{2}}).
\end{equation}
Let $\pi : \psi\inv (\lambda) \to \psi\inv (\lambda )/S^1 =
M_{\text{cut}}(\lambda , F)$ denote the orbit map.  It is easy to see
that $\pi \circ \sigma$ is onto and that $\pi \circ \sigma$ is an open
embedding on $M^{\lambda}_o$.  Moreover, since $\sigma ^* (\omega +
\omega_F) = \omega$, $\pi \circ \sigma : M^{\lambda}_o \to
M_{\text{cut}}(\lambda , F)$ is a symplectic embedding.

Similarly one checks that $\pi \circ \sigma$ induces a symplectic
embedding $j$ of $M_\lambda $ into $M_{\text{cut}}(\lambda , F)$.
Clearly, $M_{\text{cut}}(\lambda , F)\smallsetminus j(M_\lambda )
\simeq M^{\lambda}_o$. 

Assume now that $(M, \omega)$ is K\"ahler.  Then, since the symplectic
quotient of a K\"ahler manifold is K\"ahler (\cite{GS}, \cite{HKLR}),
the cut is K\"ahler as well.  We now argue that $M_{\text{cut}} (
\lambda, F) - M_\lambda$ is biholomorphic to $M^\# := \C^* \cdot
\{m\in M \mid \lambda -a < \phi(m)< \lambda\} \subset M$.

Consider the map $f:M \times \C^* \to M\times \C^*$ as in
(\ref{eqn:mapf}). We assume again as in \S~2 that the action of $S^1$
extends to an action $\tau$ of $\C^*$. Recall also that $\frac{d}{dt}
\tau_{e^t}(z) = \nabla \phi (\tau_{e^t}(z))$ for all $t\in \R$, where
$\grad \phi $ is the gradient of $\phi$ with respect to the K\"ahler
metric.  The map $f$ is biholomorphic, and intertwines the right
action (\ref{eqn:ract}) of $\C^* $ on $M\times \C^*$ with the diagonal
action (\ref{eqn:dact}) Consequently the pull-back form
$$
\tilde {\omega} : = f^* (\omega + \omega _F)
$$ 
is invariant under the action of $S^1$ on the second factor.
Moreover, this action of $S^1$ is Hamiltonian with moment map
$\tilde{\psi}$ satisfying
\begin{equation}
\tilde{\psi} = \psi_0 \circ f .
\end{equation}
where $\psi_0$ denotes the restriction of $\psi$ from $M \times \C$ to
$M \times \C^*$. Consequently $f$ induces a biholomorphic map $h$
between quotients:
$$
h : \tilde{\psi}\inv (\lambda)/S^1 \to \psi_0 \inv (\lambda)/S^1.
$$
Note that $\psi_0 \inv (\lambda)/S^1 \subset M_{\text{cut}} (\lambda,
F)$ is precisely the subset symplectomorphic to $M^{\lambda}_o$.  We
claim that $ \tilde{\psi}\inv (\lambda)/S^1$ is naturally isomorphic
to $M^\# = \C^* \cdot M^{\lambda}_o$.

Now $\tilde{\psi } \inv (\lambda) = \{ (m, z) \in M \times \C^* \mid 
\phi (\tau_{z}(m)) + H (|z|^2) = \lambda\}$. 
Since $\frac{d}{dt} \tau_{e^t}(z) = \nabla \phi (\tau_{e^t}(z))$,
$\frac{d}{dt} \phi (\tau_{e^t}(m)) = |\nabla \phi (\tau_{e^t}(m))|^2 \geq
0$.  Hence 
$$
\frac{d}{dt} \left(\phi (\tau_{e^t}(m)) + H
(|e^t|^2)\right) = |\nabla \phi (\tau_{e^t}(m))|^2 + 2 H' (e^{2t})
e^{2t} > 0.
$$
It follows that every $\C^*$ orbit intersects the level
set $\tilde{\psi}\inv (\lambda)$ in precisely one $S^1$ orbit.
Consequently
$$
\tilde{\psi} \inv (\lambda)/S^1 \simeq 
\{ m \in M \mid \phi (\tau_{e^t}(m)) + H (e^{2t}) = \lambda 
\text{ for some } t\in \R\} .
$$
We claim that the set on right hand side is $\C^* \cdot M^{\lambda}_o$.
Indeed, if for $m\in M$ there is $t\in \R$ such that $\phi (e^t \cdot
m) + H (e^{2t}) = \lambda$, then, since $H([0, \infty)) = [0, a)$ and
$H$ is strictly increasing,
$$
 \lambda -a < \phi (\tau_{e^t}(m)) < \lambda.
$$
Therefore $\tau_{e^t}(m) \in \phi \inv ((\lambda -a, \lambda)) = M^{\lambda}_o$
$\Rightarrow$ $m\in \tau_{e^{-t}}(M^{\lambda}_o) \subset \C^* \cdot
M^{\lambda}_o$.

Conversely, suppose $m\in \C^* \cdot M^{\lambda}_o$.  There are three cases to
consider: $\phi (m) \leq \lambda -a$, $\lambda -a < \phi (m) <
\lambda$ and $\lambda \leq\phi$.  Recall that 
$\lim _{t \to - \infty} H (e^{2t}) = H(0) = 0$ and that $\lim _{t \to
+ \infty} H (e^{2t}) = a$.

If $m\in \C^* \cdot M^{\lambda}_o$ and $\phi (m) \geq \lambda$ then there is
a $t\in\R$ such that $\tau_{e^t}(m) \in M^{\lambda}_o$.  So $\lim
_{t\to -\infty} \phi (\tau_{e^t}(m)) < \lambda$.  Therefore $\lim
_{t\to -\infty} (\phi (\tau_{e^t}(m))+ H(e^{2t}) < \lambda + 0$ while
$\lim _{t\to +\infty} (\phi (\tau_{e^t}(m))+ H(e^{2t})> \phi (e^0
\cdot m) + H (e^0) > \lambda$.  Consequently there is $t_1 \in \R$
with $\phi (e^{t_1} \cdot m)+ H(e^{2t_1}) = \lambda$. 

If $m\in \C^* \cdot M^{\lambda}_o$ and $\phi (m) \leq \lambda-a$ then there is
$t\in \R$ with $\tau_{e^t}(m) \in M^{\lambda}_o$.  Hence $\lim _{t\to + \infty}
(\phi (\tau_{e^t}(m)) + H (e^{2t})) > \phi (m) + a = \lambda -a + a =
\lambda$.  On the other hand
$\lim _{t\to -\infty} (\phi (\tau_{e^t}(m)) + H (e^{2t})) \leq \phi
(m) + 0 \leq \lambda -a$.  Consequently there is $t_1 \in \R$ with
$\phi (e^{t_1} \cdot m)+ H(e^{2t_1}) = \lambda$.

Finally, if $\lambda -a <\phi (m) < \lambda$ then $\lim _{t\to +
\infty} (\phi (\tau_{e^t}(m)) + H (e^{2t})) \geq \phi (m) + a > \lambda
- a + a = \lambda$ while $\lim _{t\to -\infty}(\phi (\tau_{e^t}(m))+
H(e^{2t}) < \phi (m) + 0 < \lambda$.  Consequently there is $t_1 \in
\R$ with $\phi (e^{t_1} \cdot m)+ H(e^{2t_1}) = \lambda$.
\end{proof}

\section{K\"ahler-Einstein manifolds}

Let $F$ be the function
$$
F(z) = \frac{2}{\kappa} \log(1 + \mid z \mid \sp 2), \kappa > 0.
$$
Then $\omega \sb F = \sqrt{-1}\partial {\bar \partial} F (|z|^2)$ 
is a K\"ahler-Einstein form on $\mathbb{C}$ with
structure constant $\kappa$. Hence if the symplectic form $\omega$ on $M$
is K\"ahler-Einstein with structure constant $\kappa$, so is the form
$\omega + \omega \sb F$ on $M \times \C$. Let $\psi$ be the
moment map associated to the action of $S\sp 1$ on $(M \times
\mathbb{C}, \omega + \omega \sb F)$, and let $M \sb{\cut}(\lambda, F)$
be the cut space, $\psi \sp{-1}(\lambda)/S \sp 1$. Since 
$$
a = \lim \sb{t \rightarrow +\infty} t \log(1 + t) = +\infty,
$$
$M \sb{\cut} (\lambda, F) - M \sb \lambda$ is symplectomorphic
to the open subset
$M^{\lambda}_o = \{ m \in M, -\infty < \phi(m) < \lambda\}$ of $M$ and
biholomorphic to the open subset $M \sp \# = \mathbb{C} \sp{*} \cdot
M^{\lambda}_o$ of $M$ by theorem 3.1.

The level set $Y = \psi \sp{-1}(\lambda)$ can be regarded as a principal $S
\sp 1$-bundle over the cut space:
$$
\pi: Y \to M \sb{\cut}(\lambda, F),
$$
and from the restriction of the K\"ahler metric to $Y$ one gets a
connection on this bundle and an associated curvature form $\mu$. For
$p \in M \sb{\cut}(\lambda, F)$, $\pi \sp{-1}(p)$ is an embedded circle
in $Y$, and we will denote the length of this circle measured with
respect to the K\"ahler metric by $V \sb{\eff}(p)$. For the following
we refer to \cite{BG}, $\S~11$.

\begin{thm}

If $\omega \sp \lambda$ is the K\"ahler form and $\mu \sp \lambda$ the
Ricci form on $M_{\cut}(\lambda, F)$, then
\begin{equation}
\label{eqn:ricci}
\mu \sp \lambda - 2\sqrt{-1} \partial \bar{\partial} \log V \sb{\eff}
+ c \mu = \kappa(\omega \sb \lambda + \lambda \mu)
\end{equation}
where $c$ is a constant satisfying
\begin{equation}
\label{eqn:div}
\kappa \psi = -\div \, Z + c,
\end{equation}
and $Z$ is the complex vector field on $M \times \mathbb{C}$
generating the $\mathbb{C} \sp{*}$-action. 

\end{thm}

Notice by the way that the ``$\lambda \mu$" on the right hand side of
(\ref{eqn:ricci}) has to be present for (\ref{eqn:ricci}) to be
compatible with the Duistermaat-Heckman theorem. However, the ``$c
\mu$" on the left is an artifact of (\ref{eqn:div}) which fixes the
ambiguous additive constant in the definition of the moment
map. Finally the ``$V \sb{\eff}$" on the left occurs for much the same
reason that ``effective potentials" occur for reduced Hamiltonian
systems in classical mechanics: as the contribution to the downstairs
metric of the vertical piece of the K\"ahler metric on $Y$ (see
\cite{AM}, $\S~4.5$). (Note, by the way, that the case $\kappa = 0$
can be treated similarly using $F = a|z|^2, a > 0$, as in \S~2 above.)

Here are some elementary examples of the K\"ahler cut construction.

\noindent {\bf Example 1.} Take $M = \C^n, \omega = \sqrt{-1} \sum dz_j \wedge
d\bar{z}j$, with $\tau_{e^{i\theta}}(z) = e^{i\theta}\cdot z$. Then
$\phi = |z|^2, \psi = |z|^2 + |w|^2, w \in \C,$ and the ``cut" of $M$
at $\lambda > 0$ is just $\C\mathbb{P}^n$ with $\lambda$ times the
Fubini-Study metric:
$$
\omega_\lambda = \sqrt{-1} \lambda \partial \bar{\partial}
\log(|z|^2+|w|^2).
$$

\noindent {\bf Example 2.} Take the same $M, \omega$, but with
$\tau_{e^{i\theta}}(z) = e^{-i\theta}\cdot z.$ Now $\phi = - |z|^2,
\psi = -|z|^2 + |w|^2.$ For $\lambda < 0$, $M^\lambda = \hat{\C}^n =
\C^n$ with the origin blown up. Here, in coordinates $\zeta = w\cdot
z$ on $\C^n-\{0\}$, 
\begin{equation}
\label{eqn:rfl}
\rho_\lambda = \sqrt{\lambda^2 + 4|\zeta|^2} - \lambda \log(\lambda +
\sqrt{\lambda^2 + 4|\zeta|^2}),
\end{equation}
and thus,
\begin{equation}
\label{eqn:ofl}
\omega_\lambda = \sqrt{-1} \partial \bar{\partial}
(\sqrt{\lambda^2 + 4|\zeta|^2}) - \lambda
\log(\lambda + \sqrt{\lambda^2 + 4|\zeta|^2})).
\end{equation}
Furthermore, $V_\eff = 2\pi\sqrt[4]{\lambda^2 + 4|\zeta|^2}, c = -n+1,$
and $\mu_{L,\lambda} = \sqrt{-1} \partial \bar{\partial} \log (\lambda
+ \sqrt{\lambda^2 + 4|\zeta|^2})$, so that finally
$$
Ric_\lambda = (n-1)\mu_{L,\lambda} +2 \sqrt{-1} \partial
\bar{\partial}\log V_\eff. 
$$

If we ``cut'' $\C^n$ at level $\lambda > 0$, then $M^\lambda = \C^n$,
with coordinates $\zeta = w \cdot z.$ We find the same formulas for
$\rho_\lambda$, $\omega_\lambda$ as in (\ref{eqn:rfl}),
(\ref{eqn:ofl}), though in this case the logarithmic term in
$\rho_\lambda$ extends smoothly across $\zeta = 0$. The singular case
$\lambda = 0$ has a K\"ahler metric on $\C^n -\{0\}$ with coordinates
$\zeta = w\cdot z$ as above, and $\rho_0 = 2|\zeta|$. This metric has
non-negative Ricci form which vanishes in the complex radial
directions in $\C^n$.

\noindent {\bf Example 3.} Take $M = \C\mathbb{P}^n$, with $\omega =
(n+1) \times$ the Fubini-Study K\"ahler form, so the Einstein constant
$\kappa = 1$. We will work in an affine coordinate patch $\equiv
\C^n$, with coordinates $z = (z_1,\ldots,z_n)$, so that $\omega =
\sqrt{-1} (n+1) \partial \bar{\partial} \log(1+|z|^2)$. We will blow up the
origin in the Fubini-Study metric by considering
$\tau_{e^{i\theta}}(z) = e^{-i\theta} \cdot z,$ with Hamiltonian 
\begin{equation}
\label{eqn:phifs}
\phi = (n+1) \times \frac{-|z|^2}{1+|z|^2},
\end{equation}
and cutting it with $(\C, 2 \sqrt{-1} \partial \bar{\partial}
\log(1+|w|^2))$. Thus 
\begin{equation}
\label{eqn:psifs}
\psi = -(n+1)\frac{|z|^2}{1+|z|^2} + 2\frac{|w|^2}{1+|w|^2},
\end{equation}
and so $-(n+1) < \psi < 2$ on $\C^n \times \C \ni (z,w)$. As in example 2,
for $\lambda \in (-n-1,0), M_{(\lambda, F)}, F = 2\log(1+|w|^2)$,
as in \S~3, is $\C\mathbb{P}^n$ with the origin blown up as a
complex manifold, while for $\lambda \in (0,2),$ it is just
$\C\mathbb{P}^n$ itself. The case 
$\lambda = 0$ gives a singular metric, smooth on $\C\mathbb{P}^n -
\{0\}.$ In order to begin calculating these forms, it is convenient
to do this in terms of coordinates $\zeta = w\cdot z, w$ on $\C^n -
\{0\} \times \C^*.$ From (\ref{eqn:psifs}) we solve $\psi = \lambda$ for
$|w|^2$ restricted to $\psi \inv(\lambda)$, and get
\begin{equation}
\label{eqn:wfs}
2(2-\lambda)|w|^2 = \, {\mbox{A}} + \sqrt{{\mbox{A}}^2 + {\mbox{B}}},
\end{equation}
where
$$
\begin{array}{ccl}
{\mbox{A}} & = & \lambda +(\lambda + n - 1)|\zeta|^2, \, {\mbox{and}} \\
{\mbox{B}} & = & 4(2-\lambda)(\lambda+n+1)|\zeta|^2.
\end{array}
$$
Note, curiously, that for $\lambda = -n+1$, this reduces to a
solution seen above in the case of cutting with Euclidean space. Then
\begin{equation}
\label{eqn:mufs}
\mu_{L,\, \lambda} = \sqrt{-1} \partial \bar{\partial}\log |w|^2,
\end{equation}
where $|w|^2$ is substituted using (\ref{eqn:wfs}). The qualitative
behavior of these forms at $\zeta = 0$ is the same as for the
analogous Euclidean case in example 2 above. The form
$\omega_{\lambda}$ and the effective potential are more complicated
and we do not write them down here.

\section{Cuttings by toric manifolds}
\label{sec:5.}

In this section we generalize the construction of section~\ref{sec:3.}
from $\bbC$ with a K\"ahler form defined by an $S^1$ invariant
potential to an arbitrary K\"ahler toric manifold.  The construction
is the K\"ahler counterpart of a symplectic cut with respect to a
polyhedral set (cf.\ \S~2 in \cite{LMTW}).

Let us briefly review the symplectic construction. Let $(M, \omega_M)$
be a symplectic manifold with a Hamiltonian action of a torus $G$ and
let $\varphi: M \to \fg^*$ denote an associated moment map.  Let $(X,
\omega_X)$ be a toric $G$-manifold, i.e., a (connected) symplectic
manifold with a completely integrable Hamiltonian action of the torus
$G$, {\em i.e.}, $2\dim_{\R}(G) = \dim_{\R}(M)$.  Let $\Psi: X \to
\fg^*$ denote an associated moment map.  We are mostly interested in
the cases where $X$ is a compact projective toric manifold or
$\bbC^n$.  We therefore assume that there is a convex open set $U\subset
\fg^*$ such that $\Psi (X) \subset U $ and $\Psi : X \to U$ is proper.
Then
\cite{LMTW}
\begin{enumerate}
\item The moment image $\Delta := \Psi (X)$ is locally a rational polyhedral 
set. That is for any $\eta \in \Delta$ there is a neighborhood $W$ in
$\fg^*$ and $N_1, \ldots, N_\ell \in \Z_G$ ($\Z_G$ is the integral
lattice of the torus $G$) such that 
$$
W \cap \Delta = W \cap \bigcap_{i= 1}^\ell \{ \eta \in \fg^* \mid 
\langle \eta, N_j\rangle \geq c_j\} 
$$
for some $c_1, \ldots, c_\ell \in \R$.

\item  The fibers of $\Psi$ are connected.
\end{enumerate}

\noindent Consequently, by dimension count, the fibers of $\Psi$ are
$G$-orbits. 

Assume next that (S) $\Psi : X \to \Delta$ has a continuous Lagrangian section
which is smooth over the interior $\Delta^\circ$ of $\Delta.$ This
assumption holds in a variety of contexts, for example:

\begin{enumerate}
\item if $X$ is compact, or
\item if $X$ is $\C^n$ with the standard action of the $n$-torus
$\bbT^n$ preserving the standard symplectic form, or 
\item if $X =\C$ and the symplectic form is
defined by an invariant K\"ahler potential, or 
\item if $X$ is the
symplectization of a contact toric manifold of Reeb type \cite{BG}, or
\item if $X$ is obtained from $\C^n$ by repeated symplectic cuts using
various circle subgroups of $\bbT^n$.
\end{enumerate}

The diagonal action of $G$ on $(M\times X, \omega _M - \omega _X)$ is
Hamiltonian with moment map $\Phi = \varphi -\Psi$.  We {\bf define}
the {\bf cut} of $(M,\omega_M, \varphi)$ with respect to $(X, \omega_X,
\Psi)$ to be the symplectic quotient $\bM_X := \Phi \inv (0) /G$.  If $G$ acts 
freely on $\Phi\inv (0)$, then $\bM_X$ is a smooth symplectic
manifold.  If 0 is only a regular value, then the cut space is an
orbifold.  More generally it is a symplectic stratified space in the
sense of \cite{SjL}.

Note that since $G$ is abelian, the choice of $\varphi$ and $\Psi$
involve arbitrary choices of constant vectors in $\fg^*$.  Note also
that the trivial extension of the action of $G$ on $M$ to $M\times X$
commutes with the diagonal action of $G$.  Consequently $\bM_X$ is a
Hamiltonian $G$-space.  Finally observe that the space we called $M
\sb{\cut}(\lambda, F) $ is the cut of $M$ with respect to $(\bbC,
-\sqrt{-1} \partial {\bar \partial} F (|z|^2), \lambda - H(|z|)^2)$,
the notation as in \S~3 above. 

Next we give a topological description of the cut space $\bM _X$. But
first, a bit more notation.  Since $X$ is toric the isotropy group
$G_x$ of a point $x\in X$ is connected.  Thus $G_x$ is a subtorus of
$G$.  Its Lie algebra $\fg_x$ can be read off from $\Psi (x)$ as
follows: Let $E$ be the open face of $\Delta$ containing $\Psi (x)$.
Then $\fg_x$ is the annihilator of $\Span_\R \{ E - \Psi (x)\}$.  In
particular, if $\Psi (x)$, $\Psi (x')$ are in the same face $E$ then
$G_x = G_{x'}$.  We denote this torus by $G_E$.  We now collect a 
few facts about the cut space $\bM_X$ (cf.\ \cite{Le},
\cite{LMTW}).

\begin{lemma} Let $(M, \omega_M, \varphi)$, 
$(X, \omega_X, \Psi)$, $G_E$ and  $\bM_X$ be as above.
\begin{enumerate}
\item  As a topological space 
$\bM_X: = \varphi\inv (\Delta)/\!\sim $ where $m\sim m'$
$\Leftrightarrow$ $\varphi (m) = \varphi (m')$ and $m = g\cdot m'$ for
some $g\in G_E$ where $E$ is the open face of $\Delta$ containing
$\varphi (m)$.
\item
The natural embedding of $\varphi \inv (\Delta^\circ)$ into $\bM _X$
is symplectic (as above, $\Delta^\circ$ denotes the interior of
$\Delta$).
\item
The cut space $\bM_X$ is smooth $\Leftrightarrow$ for any $m\in M$
with $\varphi (m) \in E \subset \Delta$ ($E$ is an open face) the
group $G_E$ acts freely at $m$, i.e., $G_E \cap G_m $ is trivial,
where $G_m$ denotes the isotropy group of $m$.
\end{enumerate}

\end{lemma}

\begin{proof}
Let $s :\Delta \to X$ denote a Lagrangian section of $\Psi :X \to
\Delta$, the existence of which was assumed, see (S) above.
Naturally $\Phi \inv
(0) = \{ (m, x)\in M\times X \mid \varphi (m) - \Psi (x) = 0\}$.
Hence $(m, x) \in \Phi \inv (0) $ $\Leftrightarrow$ $G\cdot x \ni 
s(\Psi (x)) = s (\varphi (m))$.  Consider the map
$$
\sigma : \varphi \inv (\Delta) \to \Phi \inv (0), 
\quad \sigma (m) = (m, s (\varphi (m))\,) .
$$ 
Its composition with the orbit map $\pi: \Phi \inv (0) \to \Phi
\inv (0)/G = \bM_X$ is surjective.  Moreover $\pi \circ \sigma$ descents to 
a homeomorphism
$$
h: \varphi \inv (\Delta) /\!\! \sim \,\,\to \bM _X .
$$
Since $s^* \omega_X = 0$,
$$
 \left( \sigma |_{\varphi \inv (\Delta^\circ)}\right)^* 
(\omega _M - \omega _X) = \omega_M |_{\varphi\inv (\Delta^\circ)} .
$$
Consequently,
\begin{equation} \label{**}
 h: \varphi\inv (\Delta^\circ) \to \bM _X
\end{equation}
is a symplectic embedding.

The cut space $\bM _X$ is smooth $\Leftrightarrow$ for all $(m, x) \in
\Phi \inv (0)$ the group $G$ acts freely at $(m, x)$.  But for $(m, x) \in
\Phi \inv (0)$ we have $\varphi (m) = \Psi (x)$, and so $G_E= G_x$.  Thus 
$G$ acts freely at $(m, x)$ $\Leftrightarrow$ $G_E \cap G_m$ is trivial.
\end{proof}

\begin{remark}
For every open face $E$ of $\Delta$ the set $M_E : = \varphi \inv (\bar{E})/\!
\sim \, \subset \bM _X$ is a submanifold fixed by $G_E$ (here
$\bar{E}$ denotes the closure of $E$ in $\Delta$), hence is
symplectic.  In particular if $E = \bar{E} = \{\mu\}$ is a vertex of
$\Delta$, then $M_E = \varphi\inv (\mu)/G$
is a symplectic quotient of $M$ at $\mu$; it is a component of the
fixed point set $(\bM _X)^G$.
\end{remark}

Let $(M, \omega _M, \varphi)$ and $(X, \omega_X, \Psi)$ be as before.
Assume further that $(M, \omega_M)$ and $(X, \omega_X)$ are K\"ahler
and that the actions of $G$ are holomorphic.  Then the cut $\bM _X$,
being a symplectic quotient, is also K\"ahler \cite{GS, HKLR}.
However, as we already observed in \S~2, the embedding $h: \varphi
\inv (\Delta^\circ) \hookrightarrow
\bM _X$ (cf.\ (\ref{**})) cannot be holomorphic relative to the complex 
structures on $M$ and $\bM_X$.  Indeed, the complement of $h (\varphi
\inv (\Delta^\circ))$ in $\bM _X$ is the union of submanifolds of the
form $M_E = \varphi \inv (\bar{E})/\!\sim $ where $E\subset \Delta$ is
a proper open face.  Each manifold $M_E$ is a component of the fixed
point set $(\bM _X)^{G_E}$, hence is $G$-invariant and K\"ahler.
Therefore $M_E$'s are preserved by the action of the complexified
group $G^\C$.  Consequently $h (\varphi\inv (\Delta^\circ)) = \bM _X
\smallsetminus \bigcup_{E\subset \Delta} M_E$ is $G^\C$ invariant. 
On the other hand there is no reason for $M^\circ := \varphi \inv
(\Delta^\circ)$ to be $G^\C$ invariant, and usually it is not.  Clearly
\begin{equation}\label{eq**}
M^\# : = G^\C \cdot M^\circ = G^\C \cdot \varphi \inv (\Delta ^\circ)
\end{equation}
is the smallest $G^\C$ invariant subset of $M$ containing $M^\circ$.
We remark:
\begin{lemma}
Let $M$, $X$, etc.\ be as above. If $G$ acts freely on $\Phi\inv (0)$
then $M^\#$ defined by (\ref{eq**}) is an open subset of $M$.
\end{lemma}

Following \cite{Sjamaar} we recall the notion of semistablility for
Hamiltonian group actions on (not necessarily integral) K\"ahler
manifolds.

\begin{definition}
Let $N$ be a K\"ahler manifold with a holomorphic Hamiltonian action
of a compact Lie group $G$ and associated moment map $\Phi: N \to
\fg^*$.  A point $x\in N$ is {\em analytically semistable} if the
closure of the $G^\bbC$ orbit through $x$ intersects the zero level
set $\Phi \inv (0)$ nontrivially ($G^\C$ denotes the complexification
of $G$).  We denote the set of semistable points in $N$ by $N_{ss}$.
\end{definition}

\begin{remark}
One can show that if the action of $G$ on $\Phi\inv (0)$ is locally
free, then $N_{ss}$ is simply $G^\bbC \cdot \Phi\inv (0)$.  In this
case one refers to the points of $N_{ss}$ as {\em stable} points and denotes
it by $N_\stable$ (cf.\ (\ref{eqn:Msharpbis} above).
\end{remark}

We will need the following property of the set of semistable points.

\begin{lemma}\label{lem-ss-ca}
Let $N$ be a K\"ahler manifold with a holomorphic Hamiltonian action
of a compact Lie group $G$ and associated moment map $\Phi: N \to
\fg^*$.  Assume that for every $x\in N$ the forward flowline of $-\nabla || \Phi||^2$,  
the negative gradient flow of the norm squared of the moment map, is
contained in a compact set.  Then the set of semistable points
$N_{ss}$ is the smallest $G^\bbC$-invariant open subset of $N$
containing $\Phi\inv (0)$.  Its complement $N \smallsetminus N_{ss}$
is a complex-analytic subset.
\end{lemma}

\begin{proof}
See \cite{Kir}, \S~4.  Compare \cite{Sjamaar}, pp. 109--110.
\end{proof}
\begin{remark}
Sjamaar ({\em op.\ cit.}) refers to moment maps $\Phi: N\to \fg^*$
with the property that for every $x\in N$ the forward flowline of $-\nabla ||
\Phi||^2$ is contained in a compact set as {\em admissible}.  All proper 
moment 
maps are admissible.  Other examples include moment maps for finite
dimensional unitary representations (Example 2.3 in \cite{Sjamaar}).
\end{remark}

We are now in a position to state and prove the main result of the
section, which is the generalization of Theorem~\ref{thm2.1} and of
Theorem~\ref{thm3.2}.

\begin{theorem}
 Let $(M, \omega_M, \varphi)$, $(X, \omega_X, \Psi)$,
 $G_E$ and $\bM_X$ be as above with $G$ acting freely on $\Phi\inv
 (0)$.  Assume further that $M$ and $X$ are K\"ahler and that $\Phi
 :M \times X \to \fg^*$ is admissible.  Then $M^\# = G^\C \cdot
 \varphi \inv (\Delta ^\circ) \subset M$ is biholomorphic to
 $h(\varphi \inv (\Delta^\circ)) \subset \bM _X$.
\end{theorem}

\begin{proof}
Fix a point $x^* \in \Psi \inv (\Delta^\circ)$. Since $X$ is connected,  the
orbit $G^\C \cdot x^* $ is all of $X^\circ : = \Psi \inv
(\Delta^\circ)$.

There are three actions of $G^\bbC$ on $M\times X$:
\begin{enumerate}

\item The trivial extension of the action of $G^\bbC$ on $M$; we denote its 
image in $\Diff (M\times X)$ by $G^\bbC _1$.
\item The trivial extension of the action of $G^\bbC$ on $X$;  we denote its 
image in $\Diff (M\times X)$ by $G^\bbC _2$.
\item The diagonal action of $G^\CC$; we denote its image in 
$\Diff (M\times X)$ by $G^\CC _d$.
\end{enumerate}

Since $G$ acts freely on $\Phi \inv (0)$,
$$
\bM _X = (M\times X)_{ss} /G^\C_d ,
$$ 
where $(M\times X)_{ss} = G^\C _d \Phi \inv (0)$ is the set of
analytically semi-stable points.  Note that by Lemma~\ref{lem-ss-ca}
$(M\times X)_{ss}$ is $G^\C_1$-invariant!

The actions of $G^\CC_1$ and $G^\CC_d$ commute.  Hence the induced
action of $G$ on $\bM_X$ is holomorphic and therefore extends to an
action of $G^\CC$. But the action of $G^\CC$ on $\bM _X$ is induced by
the action of $G^\CC_1$ on $(M\times X)_{ss}$.

Consequently $(M\times X)_{ss} $ is $G^\CC _1 \times G^\CC_2$-invariant.
Since $M\times X^\circ$ is also $G^\CC _1 \times G^\CC_2$-invariant, the
set $(M\times X^\circ)_{ss} = (M\times X^\circ) \cap (M\times X)_{ss}
$ is $G^\CC _1 \times G^\CC_2$-invariant as well.  We conclude
\begin{equation}
 (M\times X^\circ)_{ss} = \left((G^\CC_1 \times G^\CC _2) \cdot \Phi
 \inv (0)\right) \cap (M\times X^\circ)
\end{equation}

The manifold $M\times \{x^*\}$ is a cross-section for the action of
$G^\CC _d$ on $M\times X^\circ$.  Hence 
$$ 
(M\times \{x^*\} ) \cap (M\times X^\circ)_{ss} 
\simeq (M\times X^\circ)_{ss}/G^\CC _d .
$$
Thus
$$
(M\times X^\circ)_{ss} /G^\CC _d \simeq 
((G^\C_1 \times G^\CC _2) \cdot \Phi\inv (0)) \cap (M\times \{x^* \}).
$$
On the other hand,
\begin{equation*}
\begin{split}
M^{\#} & = \{ m\in M \mid \varphi (g_1\cdot m) = \Psi (x) 
\text{ for   some } g_1 \in G^\CC, x \in X^\circ\}\\
& =  \{ m\in M \mid \varphi (g_1\cdot m) -\Psi (g_2 \cdot x^*) =0 
\text{ for   some } (g_1, g_2) \in G^\CC_1\times G^\CC_2\} \\
& =  \{ m\in M \mid 
(G^\CC_1 \times G^\CC_2) \cdot (m, x^*) \cap \Phi\inv (0) \not = \emptyset\}\\
& =  \{ m\in M \mid 
 (m, x^*)\in (G^\CC_1 \times G^\CC_2) \cdot\Phi\inv (0)\}\\
& =  ((G^\CC_1 \times G^\CC_2) \cdot\Phi\inv (0)) \cap (M\times \{x^*\}).\\
\end{split}
\end{equation*}
\end{proof}

We conclude with a few very simple examples of this construction.

\vskip 2mm

\noindent {\bf {Example 1.}} One may replace, in Example 3, $\S~4$
above, $(\C,2\sqrt{-1}\partial
\bar{\partial}\log(1+|z|^2))$ by the compact toric manifold
$(\CC\mathbb{P}^1, \omega_{\mbox{FS}})$ and obtain the same cut space. 

\vskip 2mm

\noindent{\bf{Example 2.}} We note that \cite{A} has, in effect, a
large number of examples of this construction, and relates them to the
question of {\em Grauert tubes} of infinite radius. In \cite{A}, one
has a distinguished K\"ahler structure on the tangent bundle
$T\mathcal{M}$ of a real analytic Riemannian manifold $(\mathcal{M},
g)$, usually acted upon by a compact group $G$ of isometries. For
every invariant metric on $G$, one gets such a distinguished K\"ahler
strcture on $TG$ as well, and \cite{A} studies the reduction of the
product $T\mathcal{M} \times TG$, which is the same as cutting
$T\mathcal{M}$ by $TG$. Indeed, even in the case $G = S^1$, varying
the symplectic form on $TG$ by a positive scaling constant gives
interesting new examples of Riemannian manifolds $(\mathcal{M}, g')$
giving such a K\"ahler structure on $T\mathcal{M}$.

In fact, when $G = S^1$, it is easy to see that the K\"ahler
structure on $TG$ is, in the notation of this section, $X = \CC^*$
with 
$$\omega_X = c \sqrt{-1} \frac{dz \wedge d\bar{z}}{|z|^2}, c > 0.$$  
In this case, cutting {\em any} $(M, \omega_M)$ by $(X, \omega_X)$
changes neither the topology nor the complex structure on $M$,
but only alters the resultant symplectic form $\omega =
\omega(c)$. This is because the moment map on $X = \CC^*$ is proper,
and the action of $G^{\CC}$ is free and proper on $X$.

\end{document}